\documentclass[11pt,fleqn]{article}

\usepackage{amsmath,amssymb,amsthm,enumerate}

\newcommand{\p}{\partial}

\newcommand{\rank}{\mathop{\rm rank}\nolimits}

\newcommand{\CV}{\mathop{\rm CV}\nolimits}
\newcommand{\CL}{\mathop{\rm CL}\nolimits}

\newcommand{\ord}{\mathop{\rm ord}\nolimits}

\newcommand{\todo}[1][\null]{\ensuremath{\clubsuit}}

\newtheorem{theorem}{Theorem}
\newtheorem{lemma}{Lemma}
\newtheorem{corollary}{Corollary}
\newtheorem{proposition}{Proposition}
{\theoremstyle{definition}

\newtheorem{note}{Note}

\begin{document}

\par\noindent {\LARGE\bf
Local Conservation Laws\\ of Second-Order Evolution Equations\par}

{\vspace{4mm}\par\noindent \it 
Roman O. POPOVYCH~$^\dag$ and Anatoly M. SAMOILENKO~$^\ddag$ 
\par\vspace{2mm}\par}

{\vspace{2mm}\par\it 
\noindent $^{\dag,\ddag}$Institute of Mathematics of NAS of Ukraine, 3 Tereshchenkivska Str., Kyiv-4, Ukraine
\par}

{\par\noindent \it $^{\dag}$~Fakult\"at f\"ur Mathematik, Universit\"at Wien, Nordbergstra{\ss}e 15, A-1090 Wien, Austria \par}

{\vspace{2mm}\par\noindent $\phantom{^{\dag,\ddag}}$\rm E-mail: \it  $^\dag$rop@imath.kiev.ua, $^\ddag$sam@imath.kiev.ua
 \par}

{\vspace{5mm}\par\noindent\hspace*{8mm}\parbox{140mm}{\small
Generalizing results by Bryant and Griffiths [{\it Duke Math.\ J.}, 1995, V.78, 531--676], 
we completely describe local conservation laws of second-order $(1+1)$-dimensional evolution equations up to contact equivalence. 
The possible dimensions of spaces of conservation laws prove to be 0, 1, 2 and infinity. 
The canonical forms of equations with respect to contact equivalence are found for 
all nonzero dimensions of spaces of conservation laws.
}\par\vspace{4mm}}



\section{Introduction}

In the prominent paper~\cite{Bryant&Griffiths1995b} on conservation laws of parabolic equations, 
Bryant and Griffiths investigated, in particular, conservation laws of second-order $(1+1)$-dimensional evolution equations 
whose right-hand sides do not depend on~$t$. 
They proved that the possible dimensions of spaces of conservation laws for such equations are 0, 1, 2 and $\infty$. 
For each of the values 1, 2 and $\infty$, 
the equations possessing  spaces of conservation laws of this dimension were described. 
In particular, it was stated that if an evolution equation $u_t=H(x,u,u_x,u_{xx})$ has three independent conservation laws 
then this equation is linearizable.

The above results from~\cite{Bryant&Griffiths1995b} can easily be extended to 
the  general class of second-order $(1+1)$-dimensional evolution equations having the form 
\begin{equation}\label{Eq2ndOrderGenEvol}
u_t=H(t,x,u,u_x,u_{xx}),
\end{equation}
where $H_{u_{xx}}\ne0$. 
Moreover, the elimination of the restriction that the right-hand sides of the equations do not depend on~$t$
leads to an extension of the set of admissible transformations and an improvement of the transformation properties of the class. 
(Namely, the class~\eqref{Eq2ndOrderGenEvol} is normalized with respect to both point and contact transformations, 
see Section~\ref{SectionOnAdmissibleContactTransformationsOfEvolutionEquations} below.)
This allows us to essentially simplify the presentation and make more concise formulations.

In contrast to~\cite{Bryant&Griffiths1995b}, this paper does not involve differential forms. 
The conventional notions of conserved vectors and conservation laws~\cite{Olver1993} are used
(see also \cite{Popovych&Ivanova2004CLsOfNDCEs,Popovych&Kunzinger&Ivanova2008,Zharinov1986}). 

In what follows the symbol~$\mathcal L$ denotes a fixed equation from class~\eqref{Eq2ndOrderGenEvol}.
By $\CL(\mathcal L)$ we denote the space of local conservation laws of an equation~$\mathcal L$. 
It can be defined as the factor-space $\CV(\mathcal L)/\CV_0(\mathcal L)$, 
where $\CV(\mathcal L)$ is the space of conserved vectors of~$\mathcal L$ 
and $\CV_0(\mathcal L)$ is the space of trivial conserved vectors of~$\mathcal L$.
$D_t$ and $D_x$ are the operators of total differentiation with respect to the variables~$t$ and~$x$, 
$D_t=\p_t+u_t\p_u+u_{tt}\p_{u_t}+u_{tx}\p_{u_x}+\cdots$, 
$D_x=\p_x+u_x\p_u+u_{tx}\p_{u_t}+u_{xx}\p_{u_x}+\cdots$. 
Subscripts of functions denote differentiation with respect to the corresponding variables. 

The results of this paper can be summed up as follows:

\begin{theorem}\label{TheoremOnCLSOf2ndOrderEvolEqs}
$\dim\CL(\mathcal L)\in\{0,1,2,\infty\}$
for any second-order $(1+1)$-dimensional evolution equation~$\mathcal L$. 
The equation~$\mathcal L$ is (locally) reduced by a contact transformation

1) to the form $u_t=D_x\hat H(t,x,u,u_x)$, where $\hat H_{u_x}\ne0$, if and only if $\dim\CL(\mathcal L)\geqslant1$;

2) to the form $u_t=D_x^2\check H(t,x,u)$, where $\check H_u\ne0$, if and only if $\dim\CL(\mathcal L)\geqslant2$;

3) to a linear equation from class~\eqref{Eq2ndOrderGenEvol} if and only if $\dim\CL(\mathcal L)=\infty$.

\noindent 
If the equation~$\mathcal L$ is quasi-linear (i.e., $H_{u_{xx}u_{xx}}=0$) then 
the contact transformation is a prolongation of a point transformation.
\end{theorem}

\section{Admissible transformations of evolution equations}%
\label{SectionOnAdmissibleContactTransformationsOfEvolutionEquations}

It is well known~\cite{Magadeev1993} that any contact transformation mapping an equation from class~\eqref{Eq2ndOrderGenEvol} to 
an equation from the same class necessarily has the form 
\begin{equation}\label{EqContactTransOfGenEvolEqs}
\tilde t=T(t), \quad 
\tilde x=X(t,x,u,u_x), \quad 
\tilde u=U(t,x,u,u_x).
\end{equation}
The functions~$T$, $X$ and $U$ have to satisfy the nondegeneracy assumptions 
\begin{equation}\label{EqNondegeneracyAssumptionForContactTransOfGenEvolEqs}
T_t\ne0, \quad 
\rank\left(\begin{array}{ccc}X_x&X_u&X_{u_x}\\U_x&U_u&U_{u_x}\end{array}\right)=2
\end{equation}
and the contact condition  
\begin{equation}\label{EqContactConditionForContactTransOfGenEvolEqs}
(U_x+U_uu_x)X_{u_x}=(X_x+X_uu_x)U_{u_x}.
\end{equation}
The transformation~\eqref{EqContactTransOfGenEvolEqs} is uniquely prolonged to the derivatives $u_x$ and $u_{xx}$ by the formulas 
$\tilde u_{\tilde x}=V(t,x,u,u_x)$ and $\tilde u_{\tilde x\tilde x}=D_xV/D_xX$, where 
\[
V=\frac{U_x+U_uu_x}{X_x+X_uu_x}\quad\mbox{or}\quad V=\frac{U_{u_x}}{X_{u_x}}
\]
if $X_x+X_uu_x\ne0$ or $X_{u_x}\ne0$, respectively. 
The right-hand side of the corresponding transformed equation is equal to 
\begin{equation}\label{EqTransRightHandSideOfGenEvolEqs}
\tilde H=\frac{U_u-X_uV}{T_t}H+\frac{U_t-X_tV}{T_t},
\end{equation}
and $(X_u,U_u)\ne(0,0)$ in view of~\eqref{EqNondegeneracyAssumptionForContactTransOfGenEvolEqs} 
and~\eqref{EqContactConditionForContactTransOfGenEvolEqs}. 

Moreover, each of the transformations of the form~\eqref{EqContactTransOfGenEvolEqs} maps class~\eqref{Eq2ndOrderGenEvol} onto itself 
and, therefore, its prolongation to the arbitrary element~$H$ 
belongs to the contact equivalence group~$G^\sim_{\rm c}$ of class~\eqref{Eq2ndOrderGenEvol}. 
(There are no other elements in~$G^\sim_{\rm c}$.)  
In other words, the equivalence group~$G^\sim_{\rm c}$ generates the whole set of contact admissible transformations 
in class~\eqref{Eq2ndOrderGenEvol}, i.e., this class is normalized with respect to contact transformations 
(see \cite{Popovych&Kunzinger&Eshraghi2006} for rigorous definitions).
We briefly formulate the results of the above consideration in the following way. 

\begin{proposition}
Class~\eqref{Eq2ndOrderGenEvol} is contact-normalized. 
The contact equivalence group~$G^\sim_{\rm c}$ of class~\eqref{Eq2ndOrderGenEvol} is formed by the transformations~\eqref{EqContactTransOfGenEvolEqs}, 
satisfying conditions~\eqref{EqNondegeneracyAssumptionForContactTransOfGenEvolEqs} and~\eqref{EqContactConditionForContactTransOfGenEvolEqs} 
and prolonged to the arbitrary element~$H$ by~\eqref{EqTransRightHandSideOfGenEvolEqs}. 
\end{proposition}

Note that class~\eqref{Eq2ndOrderGenEvol} also is point-normalized. 
Its point equivalence group~$G^\sim_{\rm p}$ consists of the transformations 
\begin{equation}\label{EqPointTransOfGenEvolEqs}
\tilde t=T(t),\quad \tilde x=X(t,x,u),\quad \tilde u=U(t,x,u),\quad 
\tilde H=\frac\Delta{T_tD_xX}H+\frac{U_tD_xX-X_tD_xU}{T_tD_xX},
\end{equation}
where $T$, $X$ and $U$ run through the corresponding sets of smooth functions satisfying the nondegeneracy assumptions 
$T_t\ne0$ and $\Delta=X_xU_u-X_uU_x\ne0$. 

There exist subclasses of class~\eqref{Eq2ndOrderGenEvol} 
whose sets of contact admissible transformations in fact are exhausted by point transformations.

\begin{proposition}
Any contact transformation between quasi-linear equations of the form~\eqref{Eq2ndOrderGenEvol} 
is a prolongation of a point transformation.
\end{proposition}

\section{Auxiliary statements on conservation laws}%
\label{SectionOnAuxiliaryStatements}

\begin{lemma}\label{LemmaOnFormOfConsVecsOf2OrderEvolEqs}
Any conservation law of a second-order $(1+1)$-dimensional evolution equation~$\mathcal L$ contains a conserved vector~$(F,G)$
with the components $F=F(t,x,u,u_x)$ and $G=-F_{u_x}H+G^1$, where $G^1=G^1(t,x,u,u_x)$.
\end{lemma}

\begin{proof}
Let $(F,G)\in\CV(\mathcal L)$ and $\ord(F,G)=r$.
In view of the equation~$\mathcal L$ and its differential consequences, 
up to the equivalence of conserved vectors 
we can assume that $F$ and~$G$ depend only on $t$, $x$ and $u_k=\partial^k u/\partial x^k$, $k=0,\dots,r',$ where $r'\le 2r$. 
Suppose that $r'>2$.
We~expand the total derivatives in the defining relation $(D_tF+D_xG)|_{\mathcal L}=0$ for conserved vectors 
and take into account differential consequences of~$\mathcal L$ 
having the form $u_{tj}=D_x^jH$, where $u_{tj}=\partial^{j+1} u/\partial t\partial x^k$, $j=0,\dots,r'$.
Then we split the obtained condition 
\begin{equation}\label{EqDefCondForCVsOf2OrderEvolEqs}
F_t+F_{u_j}D_x^{j}H+G_x+G_{u_j}u_{j+1}=0
\end{equation}
with respect to the highest derivatives appearing in it. 
(Here the summation convention over repeated indices is used.)
Thus, the coefficients of $u_{r'+2}$ and $u_{r'+1}$ give the equations
$F_{u_{r'}}=0$ and $G_{u_{r'}}+H_{u_2}F_{u_{r'-1}}=0$ implying
\[
F=\hat F, \quad G=-S\hat F_{u_{r'-1}}u_{r'}+\hat G,
\]
where $\hat F$ and $\hat G$ are functions of $t$, $x$, $u$, $u_1$, \ldots, $u_{r'-1}$.
After selecting the terms containing $u_{r'}^2$, we additionally obtain $\hat F_{u_{r'-1}u_{r'-1}}=0$, 
i.e., $\hat F =\check F^1u_{r'-1}+\check F^0,$
where $\check F^1$ and $\check F^0$ depend at most on $t$, $x$, $u$, $u_1$,~\ldots, $u_{r'-2}$.
Consider the conserved vector with the density $\tilde F=F-D_x\Phi$ and the flux $\tilde G=G+D_t\Phi$,
where $\Phi=\int \check F^1du_{r'-2}$. It is equivalent to the initial one, and
\[
\tilde F=\tilde F(t,x,u,u_1,\ldots,u_{r'-2}), \quad
\tilde G=\tilde G(t,x,u,u_1,\ldots,u_{r'-1}).
\]

Iterating the above procedure the necessary number of times results in a conserved vector equivalent to $(F,G)$ and 
depending only on $t$, $x$, $u$, $u_1$ and $u_2$. 
Therefore, we can assume at once that $r'\le 2$.
Then the coefficients of $u_4$ and $u_3$ in~\eqref{EqDefCondForCVsOf2OrderEvolEqs} give the equations
$F_{u_2}=0$ and $G_{u_2}+H_{u_2}F_u=0$ which imply the claim.
\end{proof}

\begin{note}
Similar results are known for arbitrary $(1+1)$-dimensional evolution equations of even order~\cite{Ibragimov1985}.
In particular, any conservation law of such an equation of order~$r=2\bar r$, $\bar r\in\mathbb{N}$,  
contains the conserved vector $(F,G)$, 
where $F$ and $G$ depend only on~$t$, $x$ and derivatives of~$u$ with respect to~$x$, and
the maximal order of derivatives in~$F$ is not greater than $\bar r$.
In the proof of Lemma~\ref{LemmaOnFormOfConsVecsOf2OrderEvolEqs} 
we delibaretely used the direct method based on the definition of conserved vectors
to demonstrate its effectiveness in quite general cases. 
This proof can easily be extended to other classes of $(1+1)$-dimensional evolution equations
of even orders and some systems related to evolution equations~\cite{Popovych&Ivanova2004CLsOfNDCEs}.
\end{note}

\begin{corollary}
Any nonzero conservation law of~$\mathcal L$ is of order~1.
\end{corollary}

\begin{proof}
In view of of Lemma~\ref{LemmaOnFormOfConsVecsOf2OrderEvolEqs}, 
any conservation law of~$\mathcal L$ contains a conserved vector~$(F,G)$ 
with the components $F=F(t,x,u,u_x)$ and $G=-F_{u_x}u_t+G^1$, where $G^1=G^1(t,x,u,u_x)$. 
$(F_{u_x},G^1_{u_x})\ne(0,0)$ since otherwise condition~\eqref{EqDefCondForCVsOf2OrderEvolEqs} would imply 
that $F_u=G^1_u=0$ and, therefore, $(F,G)$ would be a trivial conserved vector. 
All trivial conserved vectors belong to the zero conservation law. 
\end{proof}

Below we consider only conserved vectors in the \emph{reduced form} which appears in Lemma~\ref{LemmaOnFormOfConsVecsOf2OrderEvolEqs}. 
For such conserved vectors, condition~\eqref{EqDefCondForCVsOf2OrderEvolEqs} is specified and expanded to  
\begin{equation}\label{EqReducedCondForCVsOf2OrderEvolEqs}
H(F_u-F_{xu_x}-F_{uu_x}u_x-F_{u_xu_x}u_{xx})+F_t+G^1_x+G^1_uu_x+G^1_{u_x}u_{xx}=0.
\end{equation}

\begin{note}\label{NoteOnTrivialConservedVectorsOf2OrderEvolEqs}
A conserved vector in reduced form is trivial if and only if its components depend at most on~$t$ and~$x$. 
If one of the components of a conserved vector in reduced form depends at most on~$t$ and~$x$ 
then the same is true for the other component.
\end{note}

\begin{lemma}\label{LemmaOnFormOfConsVecsOf2OrderEvolEqs2}
Suppose that an equation from the class~\eqref{Eq2ndOrderGenEvol} possesses a nontrivial conserved vector $(F,G)$ 
in reduced form, where additionally $F_{u_xu_x}=0$. 
Then the conserved vector $(F,G)$ is equivalent to a conserved vector~$(\tilde F,\tilde G)$ 
with $\tilde F=\tilde F(t,x,u)$ and $\tilde G=\tilde G(t,x,u,u_x)$, where $\tilde F_u\ne 0$. 
Moreover, in this case we have $H_{u_{xx}u_{xx}}=0$.
\end{lemma}

\begin{proof}
By assumption, $F=F^1u_x+F^0$ and $G=-F^1H+G^1$, where $F^1=F^1(t,x,u)$, $F^0=F^0(t,x,u)$ and $G^1=G^1(t,x,u,u_x)$. 
We put $\tilde F=F-D_x\Phi$ and $\tilde G=G+D_t\Phi$, where $\Phi=\int \check F^1du$. 
Then $\tilde F_{u_x}=0$, $\tilde G_{u_{xx}}=0$ and $(\tilde F,\tilde G)$ is a conserved vector equivalent to $(F,G)$. 
$\tilde F_u\ne 0$ since otherwise the conserved vector $(\tilde F,\tilde G)$ is trivial 
(see Note~\ref{NoteOnTrivialConservedVectorsOf2OrderEvolEqs}).
Substituting $(\tilde F,\tilde G)$ into condition~\eqref{EqReducedCondForCVsOf2OrderEvolEqs} and 
solving it with respect to~$H$, 
we obtain a linear function of~$u_x$ whose coefficients depend on~$t$, $x$ and~$u$. 
\end{proof}

\begin{corollary}\label{CorollaryOnFormOfConsVecsOf2OrderQuasiLinEvolEqs}
Any conservation law of an equation~$\mathcal L$ of the form~\eqref{Eq2ndOrderGenEvol}, where $H_{u_{xx}u_{xx}}=0$, 
contains a conserved vector $(F,G)$ with $F=F(t,x,u)$ and $G=G(t,x,u,u_x)$.
\end{corollary}

\begin{proof}
The conditions~\eqref{EqReducedCondForCVsOf2OrderEvolEqs} and $H_{u_{xx}u_{xx}}=0$ imply that 
the density of any conserved vector of~$\mathcal L$ in reduced form is linear with respect to~$u_x$. 
The claim therefore follows from Lemma~\ref{LemmaOnFormOfConsVecsOf2OrderEvolEqs2}.
\end{proof}

\begin{lemma}\label{LemmaOnFormOfConsVecsOf2OrderEvolEqs3}
If an equation~$\mathcal L$ of the form~\eqref{Eq2ndOrderGenEvol} has a nonzero conservation law 
then $H$ is a fractionally linear function in $u_{xx}$.
\end{lemma}

\begin{proof}
Suppose that $H$ is not a fractionally linear function in $u_{xx}$. 
We fix any nontrivial conserved vector $(F,G)$ of~$\mathcal L$ in reduced form. 
Such a vector exists according to Lemma~\ref{LemmaOnFormOfConsVecsOf2OrderEvolEqs}.
Splitting condition~\eqref{EqReducedCondForCVsOf2OrderEvolEqs} with respect to~$u_{xx}$ gives $F_{u_xu_x}=0$. 
Then, in view of Lemma~\ref{LemmaOnFormOfConsVecsOf2OrderEvolEqs2} 
either the function $H$ is linear in $u_{xx}$ or the conserved vector $(F,G)$ is trivial. 
This contradicts our assumption.
\end{proof}

\section{Reduction of conservation laws to canonical forms}%
\label{SectionOnReductionOfConservationLawsToCanonicalForms}

Contact equivalence transformations can be used for the reduction of equations from the class~\eqref{Eq2ndOrderGenEvol}, 
which possess nonzero conservation laws, to a special form depending on the dimension of the corresponding spaces 
of conservation laws. 
In fact, this reduction is realized via a reduction of conservation laws. 

\begin{lemma}\label{LemmaOnOneConsLawOf2OrderEvolEqs}
Any pair $(\mathcal L,\mathcal F)$, where $\mathcal L$ is an equation of the form~\eqref{Eq2ndOrderGenEvol} 
and $\mathcal F$ is a nonzero conservation law of $\mathcal L$, 
is $G^\sim_{\rm c}$-equivalent to a pair $(\tilde{\mathcal L},\tilde{\mathcal F})$, 
where $\tilde{\mathcal L}$ is an equation of the same form 
and $\tilde{\mathcal F}$ is a conservation law of $\tilde{\mathcal L}$ with characteristic~1.
\end{lemma}

\begin{proof}
Suppose that an equation~$\mathcal L$ from class~\eqref{Eq2ndOrderGenEvol} has a nonzero conservation law~$\tilde{\mathcal F}$. 
Any transformation $\mathcal T$ from $G^\sim_{\rm c}$ maps $\mathcal L$ to an equation~$\tilde{\mathcal L}$ 
from the same class~\eqref{Eq2ndOrderGenEvol} and induces a mapping from $\CL(\mathcal L)$ to $\CL(\tilde{\mathcal L})$. 
Conserved vectors of~$\mathcal L$ are transformed to conserved vectors of~$\tilde{\mathcal L}$ by the formula 
\cite{Popovych&Ivanova2004CLsOfNDCEs,Popovych&Kunzinger&Ivanova2008}
\[
\tilde F=\frac F{D_xX}, \quad \tilde G=\frac G{T_t}+\frac {D_tX}{D_xX}\frac F{T_t}.
\] 
We fix a nonzero conservation law $\mathcal F$ of~$\mathcal L$ and a conserved vector $(F,G)$ in reduced form, belonging to~$\mathcal F$, 
and immediately set $T=t$. 
The components of the corresponding conserved vector $(\tilde F,\tilde G)$ of the transformed equation~$\tilde{\mathcal L}$ 
necessarily depend at most on~$\tilde t$, $\tilde x$, $\tilde u$ and~$\tilde u_{\tilde x}$.
The conserved vector $(\tilde F,\tilde G)$ is associated with the characteristic~1 if and only if 
there exists a function $\tilde\Phi=\tilde\Phi(\tilde t,\tilde x,\tilde u,\tilde u_{\tilde x})$ such that 
$\tilde F=\tilde u+D_{\tilde x}\tilde\Phi$, i.e., in the old coordinates $D_x\Phi+UD_xX=F$, 
where  $\tilde\Phi(\tilde t,\tilde x,\tilde u,\tilde u_{\tilde x})=\Phi(t,x,u,u_x)$. 
After splitting the last equation with respect to~$u_{xx}$, we obtain the system 
\begin{equation}\label{EqCLsOfEvolEqsSystemForReductionOf1CL}
\Phi_x+UX_x+(\Phi_u+UX_u)u_x=F, \quad \Phi_{u_x}+UX_{u_x}=0.
\end{equation}
This system 
supplemented with the contact condition~\eqref{EqContactConditionForContactTransOfGenEvolEqs}
possesses the differential consequence $\Phi_u+UX_u=F_{u_x}$. 
To derive it, we need to act on the first and second equations of~\eqref{EqCLsOfEvolEqsSystemForReductionOf1CL} 
by the operators~$\p_{u_x}$ and~$\p_x+u_x\p_u$, respectively, and extract the second consequence from the first one, 
taking into account the contact condition~\eqref{EqContactConditionForContactTransOfGenEvolEqs}. 
Then system~\eqref{EqCLsOfEvolEqsSystemForReductionOf1CL} also implies  the equation $\Phi_x+UX_x=F-u_xF_{u_x}$. 
As a result, we have the system
\begin{equation}\label{EqCLsOfEvolEqsSystemForReductionOf1CL2}
\Phi_x+UX_x=F-u_xF_{u_x}, \quad \Phi_u+UX_u=F_{u_x}, \quad \Phi_{u_x}+UX_{u_x}=0.
\end{equation}
Reversing these steps shows that system~\eqref{EqCLsOfEvolEqsSystemForReductionOf1CL2}
implies~\eqref{EqContactConditionForContactTransOfGenEvolEqs} and~\eqref{EqCLsOfEvolEqsSystemForReductionOf1CL}.  
Therefore, the combined system of~\eqref{EqContactConditionForContactTransOfGenEvolEqs} and~\eqref{EqCLsOfEvolEqsSystemForReductionOf1CL} 
is equivalent to system~\eqref{EqCLsOfEvolEqsSystemForReductionOf1CL2}.  

To complete the proof, it is enough to check that for any function $F=F(t,x,u,u_x)$ with $(F_u,F_{u_x})\ne(0,0)$ 
system~\eqref{EqCLsOfEvolEqsSystemForReductionOf1CL2} has a solution $(X,U,\Phi)$ additionally satisfying 
the second condition from~\eqref{EqNondegeneracyAssumptionForContactTransOfGenEvolEqs}.

At first we consider the case $F_{u_xu_x}\ne0$ and look for solutions with $X_{u_x}\ne0$. 
The third equation of~\eqref{EqCLsOfEvolEqsSystemForReductionOf1CL2} implies that $\Phi_{u_x}\ne0$ and $U=-\Phi_{u_x}/X_{u_x}$. 
Then the two first equations take the form 
\begin{equation}\label{EqCLsOfEvolEqsSystemForReductionOf1CL3}
\Phi_x-\frac{X_x}{X_{u_x}}\Phi_{u_x}=F-u_xF_{u_x}, \quad \Phi_u-\frac{X_u}{X_{u_x}}\Phi_{u_x}=F_{u_x}.
\end{equation}
The compatibility condition of~\eqref{EqCLsOfEvolEqsSystemForReductionOf1CL3} as an overdetermined system with respect to~$\Phi$ 
is the equation 
\[
u_xF_{u_xu_x}X_x+F_{u_xu_x}X_u+(F_x-u_xF_{xu_x}-F_{uu_x})X_{u_x}=0
\]
with respect to~$X$. 
Since $F_{u_xu_x}\ne0$, this equation has a solution $X^0$ with $X^0_{u_x}\ne0$. 
The substitution of $X^0$ into~\eqref{EqCLsOfEvolEqsSystemForReductionOf1CL3} results in a compatible system with respect to~$\Phi$. 
We take a solution $\Phi^0$ of this system and put $U^0=-\Phi^0_{u_x}/X^0_{u_x}$. 
The chosen tuple $(X^0,U^0,\Phi^0)$ satisfies system~\eqref{EqCLsOfEvolEqsSystemForReductionOf1CL2}. 
The nondegeneracy condition~\eqref{EqNondegeneracyAssumptionForContactTransOfGenEvolEqs} is also satisfied. 
Indeed, suppose this was not the case. 
Then $U=\Psi(t,X)$ for some function~$\Psi$ of two arguments and system~\eqref{EqCLsOfEvolEqsSystemForReductionOf1CL} implies the equality 
\[
F=\Phi_x+\Psi X_x+(\Phi_u+\Psi X_u)u_x+(\Phi_{u_x}+\Psi X_{u_x})u_{xx}=D_x(\Phi+\textstyle\int\!\Psi\,dX),
\] 
i.e., $(F,G)$ is a trivial conserved vector. 
This contradicts the initial assumption on $(F,G)$.

If $F_{u_xu_x}=0$, in view of Lemma~\ref{LemmaOnFormOfConsVecsOf2OrderEvolEqs2} we can assume without loss of generality that $F_{u_x}=0$.
Then $F_u\ne0$. (Otherwise $(F,G)$ is a trivial conserved vector, see Note~\ref{NoteOnTrivialConservedVectorsOf2OrderEvolEqs}.)
It is obvious that the tuple $(X,U,\Phi)=(x,F,0)$ satisfies~\eqref{EqCLsOfEvolEqsSystemForReductionOf1CL2} 
and the second condition from~\eqref{EqNondegeneracyAssumptionForContactTransOfGenEvolEqs}.  
\end{proof}

\begin{corollary}\label{CorollaryOnOneConsLawOf2OrderQuasiLinEvolEqs}
Any pair $(\mathcal L,\mathcal F)$, where $\mathcal L$ is a quasi-linear equation of the form~\eqref{Eq2ndOrderGenEvol} 
and $\mathcal F$ is a nonzero conservation law of $\mathcal L$, 
is $G^\sim_{\rm p}$-equivalent to a pair $(\tilde{\mathcal L},\tilde{\mathcal F})$, 
where $\tilde{\mathcal L}$ also is a quasi-linear equation of form~\eqref{Eq2ndOrderGenEvol}
and $\tilde{\mathcal F}$ is a conservation law of $\tilde{\mathcal L}$ with characteristic~1.
\end{corollary}

\begin{proof}
In view of Corollary~\ref{CorollaryOnFormOfConsVecsOf2OrderQuasiLinEvolEqs},  
any conservation law of a quasi-linear equation of the form~\eqref{Eq2ndOrderGenEvol} possesses a conserved vector~$(F,G)$
with $F=F(t,x,u)$. 
Then the result follows from the proof of Lemma~\ref{LemmaOnOneConsLawOf2OrderEvolEqs} for the case $F_{u_x}=0$.
\end{proof}

\begin{corollary}\label{CorollaryOnFormOf2OrderEvolEqsWithOneConsLaw}
$\dim\CL(\mathcal L)\geqslant1$ if and only if the equation~$\mathcal L$ is (locally) reduced by a contact transformation
to the form $u_t=D_x\hat H(t,x,u,u_x)$, where $\hat H_{u_x}\ne0$.
The equation~$\mathcal L$ is quasi-linear if and only if
the contact transformation is a prolongation of a point transformation.
\end{corollary}

\begin{proof}
Suppose that $\dim\CL(\mathcal L)\geqslant1$. 
We fix a nonzero conservation law~$\mathcal F$ of~$\mathcal L$.  
In view of Lemma~\ref{LemmaOnOneConsLawOf2OrderEvolEqs}
the pair $(\mathcal L,\mathcal F)$ is reduced by a contact transformation~$\mathcal T$ to 
a pair $(\tilde{\mathcal L},\tilde{\mathcal F})$, 
where the equation~$\tilde{\mathcal L}$ has the form  
$\tilde u_{\tilde t}=\tilde H(\tilde t,\tilde x,\tilde u,\tilde u_{\tilde x},\tilde u_{\tilde x\tilde x})$
and $\tilde{\mathcal F}$ is its conservation law with the characteristic~1. 
If the equation~$\mathcal L$ is quasi-linear then the transformation~$\mathcal T$ is a prolongation 
of a point transformation (see Corollary~\ref{CorollaryOnFormOf2OrderEvolEqsWithOneConsLaw}).
That $\tilde{\mathcal F}$ has characteristic~1 means that the equality
$D_{\tilde t}\tilde F+D_{\tilde x}\tilde G=\tilde u_{\tilde t}-\tilde H$ is satisfied
for a conserved vector $(\tilde F,\tilde G)$ from $\tilde{\mathcal F}$. 
Therefore, up to a summand being a null divergence we have $\tilde F=\tilde u$ and $\tilde H=-D_{\tilde x}\tilde G$. 
To complete the proof, it is sufficient to put $\hat H=-\tilde G$. 

Conversely, let the equation~$\mathcal L$ be (locally) reduced by a contact transformation~$\mathcal T$
to the equation $\tilde u_{\tilde t}=D_{\tilde x}\hat H(\tilde t,\tilde x,\tilde u,\tilde u_{\tilde x})$, where $\hat H_{\tilde u_{\tilde x}}\ne0$. 
The transformed equation $\tilde u_{\tilde t}=D_{\tilde x}\hat H$ has at least one nonzero conservation law. 
This is the conservation law~$\tilde{\mathcal F}$ possessing the characteristic~1. 
The preimage of~$\tilde{\mathcal F}$ with respect to~$\mathcal T$ is a nonzero conservation law of~$\mathcal L$, 
i.e., $\dim\CL(\mathcal L)\geqslant1$. 
If $\mathcal T$ is a point transformation then the equation~$\mathcal L$ has to be quasi-linear 
as the preimage of the quasi-linear equation $\tilde u_{\tilde t}=D_{\tilde x}\hat H$ with respect to this transformation.
\end{proof}

\begin{note}\label{NoteOnConservedVectorsOfEvolEqsWithChar1}
Any conservation law of the equation~$u_t=D_x\hat H(t,x,u,u_x)$ contains a conserved vector $(F,G)$, 
where $F=F(t,x,u)$ and $G=-F_u\hat H+G^0$ with $G^0=G^0(t,x,u)$.
In this case condition~\eqref{EqReducedCondForCVsOf2OrderEvolEqs} takes the form 
$
F_t-(F_{xu}+F_{uu}u_x)\hat H+G^0_x+G^0_uu_x=0.
$ 

In particular, if additionally $F_{xu}=F_{uu}=0$ then condition~\eqref{EqReducedCondForCVsOf2OrderEvolEqs} implies 
the equations $G^0_u=0$ and $F_t+G^0_x=0$ and, therefore, $F_{tu}=0$. 
As a result, we have $F=cu+F^0(t,x)$ for some constant~$c$ and some function $F^0=F^0(t,x)$. 
This means that the conserved vector $(F,G)$ under the additional restrictions belongs to a conservation law 
which is linearly dependent with the conservation law possessing the characteristic~1.

Due to the above consideration, we can conclude that 
the space of conservation laws of the equation~$u_t=D_x\hat H(t,x,u,u_x)$ is one-dimensional if 
the right-hand side~$\hat H$ is not a fractionally linear function in $u_x$.
\end{note}

\begin{lemma}\label{LemmaOnTwoConsLawOf2OrderEvolEqs}
Any triple $(\mathcal L,\mathcal F^1,\mathcal F^2)$, where $\mathcal L$ is an equation of the form~\eqref{Eq2ndOrderGenEvol} 
and $\mathcal F^1$ and $\mathcal F^2$ are linearly independent conservation laws of~$\mathcal L$, 
is $G^\sim_{\rm c}$-equivalent to a triple $(\tilde{\mathcal L},\tilde{\mathcal F^1},\tilde{\mathcal F^2})$, 
where $\tilde{\mathcal L}$ is an equation of the same form 
and $\tilde{\mathcal F^1}$ and $\tilde{\mathcal F^2}$ are conservation laws of~$\tilde{\mathcal L}$ with the characteristics~1 and $\tilde x$.
\end{lemma}

\begin{proof}
Let the equation $\mathcal L$ possess two linearly independent conservation laws~$\mathcal F^1$ and~$\mathcal F^2$. 
We fix a conserved vector $(F^1,G^1)$ in reduced form, belonging to $\mathcal F^1$.
In view of Lemma~\ref{LemmaOnOneConsLawOf2OrderEvolEqs}, up to $G^\sim_{\rm c}$-equivalence we can assume that $F^1=u$. 
Then Lemma~\ref{LemmaOnFormOfConsVecsOf2OrderEvolEqs2} implies that $H_{u_{xx}u_{xx}}=0$ and, therefore, 
the conservation law~$\mathcal F^2$ contains a conserved vector $(F^2,G^2)$ with $F^2=F^2(t,x,u)$ and $G^2=G^2(t,x,u,u_x)$ 

We will show that there exists a point equivalence transformation of the form~\eqref{EqPointTransOfGenEvolEqs} with $T(t)=t$ such that 
the images $(\tilde F^1,\tilde G^1)$ and $(\tilde F^2,\tilde G^2)$ of the conserved vectors $(F^1,G^1)$ and $(F^2,G^2)$ 
are equivalent to the conserved vectors whose densities coincide with $\tilde u$ and $\tilde x\tilde u$, respectively. 
In other words, the conserved vectors should be transformed in such a way that 
$\tilde F^1=\tilde u+D_{\tilde x}\Phi$ and $\tilde F^2=\tilde x\tilde u+D_{\tilde x}\Psi$ 
for some functions $\Phi=\Phi(t,x,u)$ and $\Psi=\Psi(t,x,u)$. 
In the old coordinates, the conditions on $\tilde F^1$ and $\tilde F^2$ take the form $D_x\Phi+UD_xX=u$ and $D_x\Psi+XUD_xX=F^2$ 
and are split with respect to~$u_x$ to the systems  
\[
\begin{array}{l}\Phi_x+UX_x=u,\\[1ex] \Phi_u+UX_u=0\end{array}
\quad\mbox{and}\quad 
\begin{array}{l}\Psi_x+XUX_x=F^2,\\[1ex]\Psi_u+XUX_u=0.\end{array}
\]
After excluding $\Phi$ and $\Psi$ from these systems by cross differentiation, 
we derive the conditions $X_xU_u-X_uU_x=1$ and $X=F^2_u$. 
$(F^2_{xu},F^2_{uu})\ne(0,0)$ since otherwise the conservation laws~$\mathcal F^1$ and~$\mathcal F^2$ would be linearly dependent 
(see Note~\ref{NoteOnConservedVectorsOfEvolEqsWithChar1}). 
Therefore, for the value $X=F^2_u$ we have $(X_x,X_u)\ne(0,0)$. This guaranties 
the existence of a function $U=U(t,x,u)$ satisfying the equation $X_xU_u-X_uU_x=1$. 
It is obvious that the chosen functions~$X$ and~$U$ are functionally independent.
For these~$X$ and~$U$ the above systems are compatible with respect to $\Phi$ and $\Psi$. 
\end{proof}

\begin{corollary}\label{CorollaryOnTwoConsLawOf2OrderQuasiLinEvolEqs}
Any triple $(\mathcal L,\mathcal F^1,\mathcal F^2)$, where $\mathcal L$ is a quasi-linear equation of form~\eqref{Eq2ndOrderGenEvol} 
and $\mathcal F^1$ and $\mathcal F^2$ are linearly independent conservation laws of~$\mathcal L$, 
is $\smash{G^\sim_{\rm p}}$-equivalent to a triple $(\tilde{\mathcal L},\tilde{\mathcal F^1},\tilde{\mathcal F^2})$, 
where $\tilde{\mathcal L}$ is a quasi-linear equation of form~\eqref{Eq2ndOrderGenEvol} 
and $\tilde{\mathcal F^1}$ and $\tilde{\mathcal F^2}$ are conservation laws of~$\tilde{\mathcal L}$ with the characteristics~1 and $\tilde x$.
\end{corollary}

\begin{proof}
If the equation~$L$ is quasi-linear,   
$G^\sim_{\rm c}$-equivalence used in the beginning of the proof of Lemma~\ref{LemmaOnTwoConsLawOf2OrderEvolEqs}
can be replaced by $G^\sim_{\rm p}$-equivalence 
(see Corollary~\ref{CorollaryOnOneConsLawOf2OrderQuasiLinEvolEqs}). 
\end{proof}

\begin{corollary}\label{CorollaryOnFormOf2OrderEvolEqsWithTwoConsLaws}
$\dim\CL(\mathcal L)\geqslant2$ if and only if the equation~$\mathcal L$ is (locally) reduced by a contact transformation
to the form $u_t=D_x^2\check H(t,x,u)$, where $\hat H_u\ne0$.
The equation~$\mathcal L$ is quasi-linear if and only if 
the contact transformation is a prolongation of a point transformation.
\end{corollary}

\begin{proof}
In view of Lemma~\ref{LemmaOnTwoConsLawOf2OrderEvolEqs}, up to contact equivalence 
we can assume that the equation~$\mathcal L$ has the conservation laws~$\mathcal F^1$ and~$\mathcal F^2$ 
possessing the characteristics~1 and~$x$, respectively. 
(Here,  contact equivalence can be replaced by  point equivalence if the equation~$L$ is quasi-linear, 
see Corollary~\ref{CorollaryOnTwoConsLawOf2OrderQuasiLinEvolEqs}.) 
Then there exist conserved vectors $(F^1,G^1)\in\mathcal F^1$ and $(F^2,G^2)\in\mathcal F^2$ such that 
\[
D_tF^1+D_xG^1=u_t-H, \quad D_tF^2+D_xG^2=x(u_t-H).
\]
Up to the equivalence of conserved vectors, generated by adding zero divergences, we have $F^1=u$ and $F^2=xu$. 
Hence $D_xG^1=-H$ and $D_xG^2=-xH$. 
Combining these equalities, we obtain that $G^1=-D_x(G^2-xG^1)$, i.e., $H=D_x^2(G^2-xG^1)$. 
As a result, we may represent the equation~$\mathcal L$ in the form $u_t=D_x^2\check H(t,x,u)$, where $\check H=G^2-xG^1$.

Conversely, let the equation~$\mathcal L$ be reduced by a contact transformation~$\mathcal T$
to the equation $\tilde u_{\tilde t}=D_{\tilde x}^2\check H(\tilde t,\tilde x,\tilde u)$, where $\hat H_{\tilde u}\ne0$. 
The transformed equation $\tilde u_{\tilde t}=D_{\tilde x}^2\check H$ has at least two linearly independent conservation laws, e.g.,  
the conservation laws possessing the characteristics~1 and~$x$, respectively. 
Their preimages under~$\mathcal T$ are linearly independent conservation laws of~$\mathcal L$, 
i.e., $\dim\CL(\mathcal L)\geqslant2$. 
If $\mathcal T$ is a point transformation then the equation~$\mathcal L$ has to be quasi-linear 
as the preimage of the quasi-linear equation $\tilde u_{\tilde t}=D_{\tilde x}^2\check H$ with respect to this transformation.
\end{proof}

\begin{lemma}\label{LemmaOnThreeAndMoreConsLawOf2OrderEvolEqs}
$\dim\CL(\mathcal L)\geqslant3$ if and only if the equation~$\mathcal L$ is (locally) reduced by a contact transformation
to a linear equation from class~\eqref{Eq2ndOrderGenEvol}. 
The equation~$\mathcal L$ is quasi-linear if and only if 
the contact transformation is a prolongation of a point transformation.
\end{lemma}

\begin{proof}
Let $\dim\CL(\mathcal L)\geqslant3$. 
In view of Corollary~\ref{CorollaryOnFormOf2OrderEvolEqsWithTwoConsLaws}, 
the equation~$\mathcal L$ can be assumed, up to $G^\sim_{\rm c}$-equivalence, 
to have the representation $u_t=D_x^2\check H(t,x,u)$, where $\hat H_u\ne0$. 
Here, $G^\sim_{\rm p}$-equivalence can be used instead of $G^\sim_{\rm c}$-equivalence if $\mathcal L$ is a quasi-linear equation.
Then condition~\eqref{EqReducedCondForCVsOf2OrderEvolEqs} implies that 
each conservation law of~$\mathcal L$ contains a conserved vector $(F,G)$, 
where $F=F(t,x,u)$ and $G=-F_u\hat H+G^0$ with $G^0=G^0(t,x,u)$ (cf. Note~\ref{NoteOnConservedVectorsOfEvolEqsWithChar1}). 
Additionally, the function~$F$ and~$G^0$ have to satisfy the equations
\[
F_{uu}=0, \quad F_u\check H_{xu}-F_{xu}\check H_u+G^0_u=0, \quad F_t+F_u\check H_{xx}+G^0_x=0.
\] 
The first equation gives that, up to the equivalence of conserved vectors, generated by adding zero divergences,
$F=fu$ with some function~$f=f(t,x)$.
Exclusion of $G^0$ from the other equations by cross differentiation leads to the condition $f_t+f_{xx}\check H_u=0$. 
If we would have $\check H_{uu}\ne0$, this condition would imply $f_t=f_{xx}=0$, i.e., $f\in\langle1,x\rangle$. 
In other words, any conservation law of~$\mathcal L$ would be a linear combination of 
the conservation laws possessing the characteristics~1 and~$x$ if $\check H_{uu}\ne0$. 
Therefore, since $\dim\CL(\mathcal L)\geqslant3$, the case $\check H_{uu}\ne0$ is impossible. 
The condition $\check H_{uu}=0$ is equivalent to the equation $u_t=D_x^2\check H(t,x,u)$ being linear. 

Conversely, suppose that the equation~$\mathcal L$ is reduced by a contact transformation~$\mathcal T$
to a linear equation~$\tilde{\mathcal L}$ from class~\eqref{Eq2ndOrderGenEvol}. 
The space of conservation laws of any linear equation (with sufficiently smooth coefficients) is infinite-dimensional.   
Therefore, the space $\CL(\mathcal L)$ is infinite-dimensional as the preimage of the infinite-dimensional space $\CL(\tilde{\mathcal L})$ 
with respect to the one-to-one mapping from $\CL(\mathcal L)$ onto $\CL(\tilde{\mathcal L})$, generated by~$\mathcal T$. 
If $\mathcal T$ is a point transformation then the equation~$\mathcal L$ has to be quasi-linear 
as the preimage of the linear equation~$\tilde{\mathcal L}$ with respect to this transformation.
\end{proof}

\section{Examples}%
\label{SectionOnExamples}

Conservation laws of different subclasses of class~\eqref{Eq2ndOrderGenEvol} were classified in a number of papers 
(see, e.g., \cite{Bryant&Griffiths1995b,Ivanova&Popovych&Sophocleous2004,Popovych&Ivanova2004CLsOfNDCEs,Popovych&Kunzinger&Ivanova2008} 
and the references therein).
All known results perfectly agree with Theorem~\ref{TheoremOnCLSOf2ndOrderEvolEqs}. 

Thus, both local and potential conservation laws of nonlinear diffusion--convection equations of the general form 
\begin{equation} \label{EqDC}
u_t=(A(u)u_x)_x+B(u)u_x,
\end{equation}
where $A=A(u)$ and $B=B(u)$ are arbitrary smooth functions of $u$ and $A(u)\ne0$, 
were exhaustively investigated in~\cite{Popovych&Ivanova2004CLsOfNDCEs}. 
The point equivalence group $G^{\sim}$ of the class~\eqref{EqDC} is formed by the transformations
\[
\tilde t=\varepsilon_4t+\varepsilon_1, \quad
\tilde x=\varepsilon_5x+\varepsilon_7 t+\varepsilon_2, \quad
\tilde u=\varepsilon_6u+\varepsilon_3, \quad
\tilde A=\varepsilon_4^{-1}\varepsilon_5^2A, \quad
\tilde B=\varepsilon_4^{-1}\varepsilon_5B-\varepsilon_7,
\]
where $\varepsilon_1,$ \dots, $\varepsilon_7$ are arbitrary constants,
$\varepsilon_4\varepsilon_5\varepsilon_6\ne0.$
Any equation from class~\eqref{EqDC} possesses the conservation law $\mathcal F^0$ 
whose density, flux and characteristic are 
\[
\mathcal F^0=\mathcal F^0(A,B)\colon\quad F=u,\quad G=-Au_x-\breve B,\quad \lambda=1.
\]
A complete list of $G^\sim$-inequivalent equations~\eqref{EqDC} having
additional (i.e., linearly independent of~$\mathcal F^0$) conservation laws
is exhausted by the following ones
\[\arraycolsep=0ex\begin{array}{lllll}
B=0,\quad&\mathcal F^1=\mathcal F^1(A)\colon\quad& F=xu,& \quad G=\breve A-xAu_x,&\quad \lambda=x;\\[1ex]
B=A,\quad&\mathcal F^2=\mathcal F^2(A)\colon\quad& F=e^xu,& \quad G=-e^xAu_x,&\quad \lambda=e^x;\\[1ex]
A=1,\quad& B=0,\quad\mathcal F^3_h\colon\quad&  F=h u,& \quad G=h_xu-h u_x,&\quad \lambda=h.
\end{array}\]
where $\breve A=\int\! A(u) du$, $\breve B=\int\! B(u) du$ and 
$h=h(t,x)$ runs through the set of solutions of the backward linear heat equation $h_t+h_{xx}=0$.
(Along with constrains for $A$ and $B$ the above table also contains complete lists of
densities, fluxes and characteristics of the additional conservation laws.)
Therefore, all possible nonzero dimensions of spaces of conservation laws of evolution equations are realized in the class~\eqref{EqDC}. 
Moreover, excluding one case, the equations listed above are already represented in the corresponding canonical forms 
which are described in Theorem~\ref{TheoremOnCLSOf2ndOrderEvolEqs}.
To reduce an equation from class~\eqref{EqDC} with $B=A$ to the canonical form of evolution equations 
possessing two linearly independent conservation laws (item~2 of Theorem~\ref{TheoremOnCLSOf2ndOrderEvolEqs}), 
according to the proof of Lemma~\ref{LemmaOnTwoConsLawOf2OrderEvolEqs}, 
we have to apply the transformation $\tilde t=t$, $\tilde x=e^x$ and $\tilde u=e^{-x}u$. 
The transformed equation $\tilde u_{\tilde t}=D_{\tilde x}^2\breve A(\tilde x\tilde u)$ 
does not belong to the class~\eqref{EqDC} but is represented in the canonical form. 

More generally, suppose that an evolution equation has two linearly independent conservation laws 
whose characteristics $\lambda^1$ and $\lambda^2$ depend at most on~$t$ and~$x$. 
Then a transformation reducing this equation to the canonical form is 
$\tilde t=t$, $\tilde x=\lambda^2/\lambda^1$ and $\tilde u=\lambda^1u/(\lambda^2/\lambda^1)_x$. 
This gives a simple way for finding the corresponding transformations, e.g.,  
in the class of variable coefficient diffusion--convection equations of the form
$f(x)u_t=(g(x)A(u)u_x)_x+h(x)B(u)u_x$. 
The local conservation laws of such equations were investigated in~\cite{Ivanova&Popovych&Sophocleous2004}.

As nontrivial examples on case~3 of Theorem~\ref{TheoremOnCLSOf2ndOrderEvolEqs}, we consider the linearizable equations 
$\mathcal L_1$: $u_t=u_x{}^{\!\!-2}u_{xx}$ and $\mathcal L_2$: $u_t=-u_{xx}{}^{\!\!-1}$. 
They are the first and second level potential equations of the remarkable diffusion equation $u_t=(u^{-2}u_x)_x$ 
and are reduced to the linear heat equation $\tilde u_{\tilde t}=\tilde u_{\tilde x\tilde x}$ 
by the (point) hodograph transformation $\tilde t=t$, $\tilde x=u$ and $\tilde u=x$ 
and the (contact) Legendre transformation $\tilde t=t$, $\tilde x=u_x$ and $\tilde u=xu_x-u$, respectively.  
The spaces $\CL(\mathcal L_1)$ and $\CL(\mathcal L_2)$ are infinite-dimensional. 
The space $\CL(\mathcal L_1)$ consists of the conservation laws
with the conserved vectors $(F,G)=(\sigma,\sigma_\omega u_x{}^{\!\!-1})$ and the characteristics $\lambda=\sigma_\omega$, where $\omega=u$. 
The space $\CL(\mathcal L_2)$ is formed by the conservation laws
with the conserved vectors $(F,G)=(\sigma,\sigma_\omega u_{xx}{}^{\!\!-1})$ and the characteristics $\lambda=\sigma_tu_{xx}$, where $\omega=u_x$. 
In both the cases the parameter-function $\sigma=\sigma(t,\omega)$ runs through the solution set 
of the backward linear heat equation $\sigma_t+\sigma_{\omega\omega}=0$.

The unified representations of equations possessing conservation laws are important for 
a successful study of the potential frame (potential systems, potential conservation laws and potential symmetries) 
for the class~\eqref{EqDC}, confer also the discussion in the next section.

\section{Conclusion}%
\label{SectionConclusion}

In this paper we have presented the classification of conservation laws of general second-order (1+1)-dimensional evolution equations. 
The classification list is very compact.  
In addition to the odd order and the evolution structure of the equations under consideration, 
the simplicity of the classification result is explained by the normalization of the class of these equations 
with respect to contact transformations. 
(The class considered in~\cite{Bryant&Griffiths1995b} is not normalized.)

The classification of local conservation laws leads to the complete description of first-level potential systems 
of evolution equations. 
The contact equivalence group $G^\sim_{\rm c}$ of the class~\eqref{Eq2ndOrderGenEvol} generates an equivalence relation 
on the corresponding set of potential systems~\cite{Popovych&Ivanova2004CLsOfNDCEs,Popovych&Kunzinger&Ivanova2008}. 
Up to this equivalence relation and the equivalence of conserved vectors, 
the first-level potential systems of those equations non-linearizable by contact transformations 
are exhausted by the systems
\begin{gather}\label{EqPotSystems1Of2ndOrderGenEvol}
v_x=u, \quad v_t=\hat H,
\\\intertext{where  $\hat H=\hat H(t,x,u,u_x)$ and $\hat H_{u_x}\ne0$, and}
\label{EqPotSystems1xOf2ndOrderGenEvol}
v^1_x=u, \quad v^1_t=D_x\check H, \quad 
v^2_x=xu, \quad v^2_t=xD_x\check H-\check H,
\end{gather}
where  $\check H=\check H(t,x,u)$ and $\check H_u\ne0$.

Each system of the form~\eqref{EqPotSystems1Of2ndOrderGenEvol}
is constructed with a single conserved vector in reduced form, associated with the characteristic~1. 
The corresponding potential equation is $ v_t=\hat H(t,x,v_x,v_{xx})$.

Each system of the form~\eqref{EqPotSystems1xOf2ndOrderGenEvol}
is constructed with a pair of conserved vectors in reduced form, associated with the characteristics~1 and~$x$. 
It can formally be represented as the second-level potential system
\begin{equation}\label{EqPotSystems1xOf2ndOrderGenEvolAsSecondLevelPotSystem}
v^1_x=u, \quad w_x=v^1, \quad w_t=\check H(t,x,u),
\end{equation}
where $w=xv^1-v^2$. 
The equation $v^1_t=D_x\check H$ is a differential consequence of the second and third equations 
of~\eqref{EqPotSystems1xOf2ndOrderGenEvolAsSecondLevelPotSystem} and can be omitted from the canonical representation.
The potential equation associated with~\eqref{EqPotSystems1xOf2ndOrderGenEvolAsSecondLevelPotSystem} is $w_t=\hat H(t,x,w_{xx})$.
In spite of formally belonging to the second level of potential systems,
the representation~\eqref{EqPotSystems1xOf2ndOrderGenEvolAsSecondLevelPotSystem} has a number of advantages 
in comparison with the representation~\eqref{EqPotSystems1xOf2ndOrderGenEvol}. 

An exhaustive study of the potential frame for linear second-order (1+1)-dimensional evolution equations, 
including potential systems, potential conservation laws, usual and generalized potential symmetries of all levels, 
was presented in~\cite{Popovych&Kunzinger&Ivanova2008}. 

The Lie symmetries of the first-level potential systems~\eqref{EqPotSystems1Of2ndOrderGenEvol} and~\eqref{EqPotSystems1xOf2ndOrderGenEvol} 
are the first-level potential symmetries of equations from the class~\eqref{Eq2ndOrderGenEvol}.
System~\eqref{EqPotSystems1xOf2ndOrderGenEvol} can be replaced by 
system~\eqref{EqPotSystems1xOf2ndOrderGenEvolAsSecondLevelPotSystem} since these systems are point-equivalent.
To investigate Lie symmetries of~\eqref{EqPotSystems1Of2ndOrderGenEvol} and~\eqref{EqPotSystems1xOf2ndOrderGenEvolAsSecondLevelPotSystem}, 
results of \cite{Magadeev1993}
(resp. \cite{Abramenko&Lagno&Samojlenko2002,Basarab-Horwath&Lahno&Zhdanov2001,Lagno&Samojlenko2002})  
on the classification of contact (resp. Lie) symmetries of equations from the class~\eqref{Eq2ndOrderGenEvol} 
with respect to the corresponding contact (resp. point) equivalence group can be used. 
The simplest case of this strategy was discussed in~\cite{Zhdanov&Lahno2005}.

\looseness=-1
The iterative application of the procedure of finding conservation laws to potential systems together 
with the subsequent construction of potential systems of the next level gives 
a description of universal Abelian coverings~\cite{Bocharov&Co1997} 
(or extensions by conservation laws in the terminology of~\cite{Bryant&Griffiths1995b}).
See also~\cite{Kunzinger&Popovych2008,Marvan2004} for a definition of Abelian coverings and 
\cite{Sergyeyev2000} for a discussion of universal Abelian coverings of evolution equations. 
As a next step we will complete the study of universal Abelian coverings for equations from the class~\eqref{Eq2ndOrderGenEvol}, 
using the equivalence relation generated by the contact equivalence group and other techniques. 
These results will form the subject of a forthcoming paper.

\subsection*{Acknowledgements}

The authors are grateful to V.~Boyko, M.~Kunzinger and A.~Sergyeyev for productive and helpful discussions.
The research of ROP was supported by START-project Y237 and Stand-alone project P20632 of the Austrian Science Fund. 
The authors also wish to thank the referee for his/her suggestions for the improvement of this paper.

\end{document}